\documentclass[11pt,reqno]{amsart}

\usepackage{amsfonts}
\usepackage{latexsym}
\usepackage{amssymb}

\newcommand{\C}{\mathbb C}

\newcommand{\Z}{\mathbb Z}

\newtheorem{thm}{Theorem}[section]
\newtheorem{prop}[thm]{Proposition}
\newtheorem{lem}[thm]{Lemma}
\newtheorem{rem}[thm]{Remark}
\newtheorem{cor}[thm]{Corollary}

\newtheorem{defn}[thm]{Definition}

\begin{document}

\title[]{Degeneration of Schubert Varieties of 
$SL_{n}/B$ to Toric Varieties}

\author[]{Raika Dehy 
\and
Rupert W.T. Yu}
\address{D\'epartement de Math\'ematiques\\
Universit\'e de Cergy-Pontoise \\
2 avenue Adolphe Chauvin \\
95032 Cergy, France}
\email{dehy@math.pst.u-cergy.fr}
\address{UMR 6086 CNRS\\
D\'epartement de Math\'ematiques\\
Universit\'e de Poitiers\\
Boulevard Marie et Pierre Curie\\
T\'el\'eport 2 - BP 30179\\
86962 Futuroscope Chasseneuil cedex, France}
\email{yuyu@mathlabo.univ-poitiers.fr}

\begin{abstract}
Using the polytopes defined in an earlier paper, we show in this paper
the existence of degeneration of a large class of Schubert varieties 
of $SL_{n}$ to toric varieties by extending the method used by 
Goncuilea and Lakshmibai for a miniscule $G/P$ to Schubert varieties 
in $SL_{n}$.
\end{abstract}


\maketitle

\section*{Introduction}
\setcounter{equation}{0}

In this paper, we complete our programme stated in \cite{kn:D-Y} to prove the 
existence of degenerations of certain Schubert varieties of $SL_{n}$ into toric 
varieties, thus generalizing the results of Goncuilea and Lakshmibai
\cite{kn:G-L}. 
For example, we are able to settle all the Schubert varieties in $SL_{3}$ here.

The essential idea is that we use the polytopes defined in \cite{kn:D-Y} to 
construct a distributive lattice, and extend the method used
by Goncuilea and Lakshmibai \cite{kn:G-L} for miniscule $G/P$ to Schubert 
varieties in $SL_{n}$. Although they also prove the existence of 
degenerations for $SL_{n}/P$ (and also Kempf varieties) in the same
paper, their approach is different from the one for a
miniscule $G/P$.

Since all the ingredients used here are based on standard monomials,
we expect that it can be adapted in the other types.  
However, the difficult part is to construct a suitable distributive lattice 
and we shall make it more precise below.

Let $G=SL_{n+1}$, $B$ be a Borel subgroup and $W$ be the Weyl group of $G$
which is the symmetric group of $n+1$ letters. Let $\alpha_{i}$, 
$i=1,\cdots ,n$, be the corresponding set of simple roots so that
$\langle \alpha_{i},\alpha_{j}^{\vee}\rangle =a_{i,j}$ 
where $(a_{i,j})_{i,j}$ is
the Cartan matrix, $s_{i}$ the corresponding simple reflections
in $W$ and let $\omega_{i}$ be the corresponding fundamental weights. 
Denote also by $\ell (-)$ and $\preceq$ the length
function and the Bruhat order on $W$. 

Recall that for  $w \in W$,  the Demazure module
$E_w (\lambda)$ is the $\mathfrak{b}$-module
$U(\mathfrak{b}) v_{w \lambda}$, where $\mathfrak{b}$
is the Lie algebra of $B$, $U(\mathfrak{b})$ its envelopping
algebra and $v_{w \lambda}$ a vector of extremal 
weight $w \lambda$ of the irreducible representation
$V(\lambda)$ of highest weight $\lambda =\sum_{i=1}^{n}k_{i}\omega_{i}$,
$k_{i}\geq 0$. 
Under certain conditions on $w$, in \cite{kn:D-Y}, we constructed $n$ 
polytopes  $\Delta_1, \ldots , \Delta_n$, where $n$ is the
rank of $G$, such that the number of lattice points in the Minkowski 
sum $\sum_{i=1}^n k_i \Delta_i = \{  \sum_{i=1}^n \sum_{j=1}^{k_{i}}
x_{ij} \mid x_{ij} \in \Delta_i \}$ is equal to
the dimension of $E_w (\lambda)$.
The polytopes $\Delta_1, \ldots , \Delta_n$ define a toric 
variety $X$ equipped with $n$ line bundles 
$\mathcal{L}_i, i=1 , \ldots , n$ 
(see \cite{kn:T}). The aim of this paper is
to degenerate the Schubert variety $S(w) =\overline{BwB}/B$ equipped
with line bundles 
$\mathcal{L}_{\omega_i}= \overline{BwB} \times_{B} \C_{\omega_i}$ into $X$ 
equipped with $\mathcal{L}_i$.

We consider the homogeneous coordinate ring of a multicone
over $S_w$. This multicone is the
$B \tau B$-orbit of $\bigoplus_{i=1}^n \C v_{\omega_i}$ in
$\bigoplus_{i=1}^n V (\omega_i)$, and its coordinate ring is
$R = \bigoplus_{\lambda\, {\rm dominant}}
H^0 (S_w , \mathcal{L}_{\lambda})$, where
$\mathcal{L}_{\lambda} = \bigotimes_{i=1}^n
\mathcal{L}_{\omega_i}^{\otimes k_i}$ with 
$\lambda = \sum_i k_i \omega_i$. In 
\cite{kn:K-R}, it has been shown that the map
$$
\bigoplus_{k_1, \ldots , k_n \geq 0} 
\bigotimes_{i=1}^n {\rm Sym}^{k_i} H^0 (S_w , \mathcal{L}_{\omega_i})
\rightarrow R
$$ 
is surjective and its kernel $I$
is a multigraded ideal generated by elements of
degree $(k_1, \ldots ,k_n)$ with $\sum_i k_i =2$.
On the other hand, we encounter an analogous situation
considering the toric variety $X$ defined by the polytopes 
$\Delta_1, \ldots , \Delta_n$. Let $\mathcal{B}_{k_1, 
\ldots , k_n}$ be the vector space over $\C$ generated by
$x^{\alpha}$, $\alpha$ a lattice point in $\sum_{i=1}^n k_i \Delta_i$. Then
$S = \bigoplus_{k_1, \ldots ,
k_n \geq 0} \mathcal{B}_{k_1, \ldots , k_n }$ is 
the homogeneous coordinate ring of a multicone 
over the toric variety $X$, and 
$\mathcal{B}_{k_1, \ldots , k_n} = H^0 (X, \bigotimes_{i=1}^n 
\mathcal{L}_{\omega_i}^{\otimes k_i})$. Moreover,
since the polytopes $\Delta_i$ can be triangulized by 
simplices of minimal volume, that is of volume $1/ 
(\dim \Delta_{i})!$, the map 
$\bigoplus_{k_1, \ldots , k_n \geq 0} \bigotimes_{i=1}^n
{\rm Sym}^{k_i} H^0 (X , \mathcal{L}_i)
\rightarrow S$ is surjective and its kernel $J$
is a multigraded ideal generated by elements of
degree $(k_1, \ldots , k_n)$ with $\sum_i k_i =2$;
in other words $S = \C [ x^{\alpha_{i,j}}]/J$
where $H := \{ \alpha_{i,j} \}$ is the set of all 
lattice points in polytopes $\Delta_i$, $i =1, \ldots , n$.

The basic idea is that one can put a structure of 
a distributive lattice  on the set $H$, of lattice points of 
$\Delta_1, \ldots , \Delta_n$. This distributive lattice, 
denoted $H$ equipped with operations $\vee$, $\wedge$, is
such that for $\alpha , \beta \in H$ we have
$\alpha + \beta = \alpha \vee \beta + \alpha  \wedge
\beta$. Hence that the algebra $\C [H]/ I(H)$,
where $I(H)$ is the homogeneous ideal generated by
$x_{\alpha} x_{\beta} = x_{\alpha \vee \beta} x_{\alpha  
\wedge \beta}$, is the ring $S$. Therefore using 
theorem~\ref{thm:degeneration} proved in \cite{kn:G-L}, 
one obtains a flat deformation of $R$ to $\C [H]/ I(H)$ which is the 
homogeneous coordinate ring of a multicone over the
toric variety $X$.

The paper is organized as follows. In section 1, we
recall results from \cite{kn:D-Y}. The theorem of
degeneration of \cite{kn:G-L} is stated in section 2. 
Sections 3,4 and 5 are devoted to showing that
the conditions of the theorem are satisfied. Finally
in section 6, we discuss briefly which Schubert varieties
fall into our context.

We shall use the above notations throughout this paper.

\section{Distributive lattice on ${\mathcal W}^{w}$}
\setcounter{equation}{0}
\label{sec:prelim}

For a fundamental weight $\omega_i$, $i =1, \ldots , n$, 
let $W_{\omega_i}$ be the subgroup of the Weyl group $W$, stabilizing 
$\omega_i$, that is $W_{\omega_i}=\{\tau\in W \mid\tau (\omega_i)=\omega_i \}$. 
Denote the quotient $W/W_{\omega_i}$ by $W_i$. 
The set $W_i$ can, on the one hand, be identified with the subset of $W$ 
consisting of elements $\tau$ such that $\tau \preceq \tau s_{\alpha_j}$ for 
$j \not = i$, {\it i.e.} the set of minimal representatives and, 
on the other hand, with the set of $i$-tuples $(r_1, \ldots, r_i)$ such that
$0 \leq r_1 < \cdots < r_i \leq n$. The connection between
these two identifications is that $(r_1, \ldots , r_i)$ 
corresponds to $s(r_1,1) s(r_2,2) \cdots s(r_i,i)$
where $s(a,b) = s_{a} s_{a-1} \cdots s_b$.
The induced Bruhat order on  $W_i$, which we shall
also denote by $\preceq$ can be expressed under the above identifications by
$\underline{a} = (a_1, \ldots , a_i)\preceq \underline{b} =(b_1, \ldots , b_i)$
if and only if $a_k \leq b_k$, $1 \leq k \leq i$. 
Furthermore, $W_i$ becomes a distributive lattice 
(for generalities on distributive lattices,
see \cite{kn:H} or section~2 of \cite{kn:G-L}) under $\preceq$ where 
\begin{equation}
\begin{array}{rl}
\label{eq:genmax-min}
&\underline{a} \vee \underline{b}
= (\max \{a_1, b_1 \}, \ldots , \max \{a_i,b_i\}) \\ 
{\rm and} &\underline{a} \wedge \underline{b} = (\min \{a_1, b_1\}, 
\ldots , \min \{a_i,b_i\}).
\end{array} 
\end{equation}

Recall (see for example~\cite{kn:D-Y}) that any $w \in W$ has a
unique factorization in the form 
$s(a_1,b_1) s(a_2,b_2) \cdots s(a_k,b_k)$ with 
$1 \leq a_1 < a_2 < \cdots <a_k \leq n$.
We shall be interested in the $w$'s satisfying
$b_1 \geq b_2 \geq \cdots \geq b_k$.

For an element $w \in W$, let 
$W_i^w =\{ \tau \in W_i \mid \tau \preceq \overline{w} \}$,
where $\overline{w}$ is the representative of $w$ in
$W_i$. Denote by ${\mathcal W}^w := \coprod_{i=1}^n W_i^w$.
Let us recall the following partial order from
section~8 of \cite{kn:D-Y}.

\begin{defn}\label{defn:bar,tilde}
Let $i \leq j$ and $w = s(a_1,b_1) \cdots s(a_k,b_k)$,
with $1 \leq a_1 < \cdots < a_k \leq n$ and $b_1 \geq b_2 \geq \cdots \geq b_k$.
For $\phi = (r_1, \ldots , r_i)\in W_i^w$, we define 
$$
\breve{\phi}:= (0, 1, \ldots , j-i-1, \breve{r}_{j-i+1}, \ldots ,
\breve{r}_j) \in W_j^w
$$ 
where $\breve{r}_k = \max \{k-1, r_{k-j+i}\}$, $j-i+1\leq k \leq j$ and for 
$\tau = (t_1, \ldots , t_j)\in W_j^w$, let 
$$
\tilde{\tau} :=(t_{j-i+1}, \ldots , t_j) \in W_i^w.  
$$
We say that $\phi \preceq_w \tau$
if $\breve{\phi} \preceq \tau$, or equivalently if $\phi\preceq \tilde{\tau}$, 
and we define
$\tau \vee \phi :=\tau \vee \breve{\phi} \in W_j^w$ and
$\tau \wedge \phi := \tilde{\tau} \wedge \phi \in W_i^w$
(see Eq.~\ref{eq:genmax-min}).
\end{defn}

A simple consequence of the definition is the following lemma.

\begin{lem}\label{distri}
Let $w$ be as in definition~\ref{defn:bar,tilde}. Then
together with the above operations, ${\mathcal W}^{w}$ is a distributive
lattice. 
\end{lem}

An essential property of this partial order is the following
theorem proved in \cite{kn:D-Y}.

\begin{thm}\label{stdness}
We have $\phi \preceq_w \tau$ in ${\mathcal W}^w$ if and only if there exist 
liftings $\phi^{\prime}$, $\tau^{\prime}$ in $W$ of 
$\phi, \tau$ such that $\phi^{\prime} \preceq \tau^{\prime}\preceq w$.
\end{thm}

As we shall see in the next sections, this is used extensively
in the proof.

\begin{rem}\label{rem:polytope}
{\rm 
In \cite{kn:D-Y}, we constructed for each fundamental weight $\omega_{i}$, 
a polytope $\Delta_i^{w}$ such that the number of lattice points in the 
Minkowski sum $\sum_{i=1}^n k_i \Delta_i^{w}$ is equal to 
$\dim E_w (\sum_{i=1}^n k_i \omega_i)$. The set of vertices of the polytope 
$\Delta_i$ is indexed by the 
set $W_i^w$ and these are the only lattice points of $\Delta_i$. Moreover 
considering  $\phi , \tau \in \mathcal{W}^w$ as vertices, we have
$\phi + \tau = \phi \vee \tau + \phi \wedge \tau$. 
The polytopes $\Delta_i$ have also the
important property that they can be triangulized by
simplices of minimal volume so that a lattice point of 
$\sum_{i=1}^n k_i \Delta_i$ can be written as the sum of $k_1$
lattice points of $\Delta_1$ and $k_2$ lattice points of 
$\Delta_2$ and so on. This property gives information on the 
generators of the toric ideal defined by the $\Delta_i$.
}
\end{rem}

We shall end this section by proving certain facts concerning
$\tau \vee \phi$ and $\tau \wedge \phi$ which will be needed 
throughout the paper. These are generalizations of certain results
obtained in \cite{kn:G-L}. Let us suppose that $w$ is as in
definition~\ref{defn:bar,tilde}. 

\begin{lem}\label{lem:equalsum1}
Let $j\geq i$ and $\phi \in W_i^w$, $\tau \in W_j^w $ be two non-comparable
elements in ${\mathcal W}^w$. Let $\sigma = \tau \vee \phi$
and $\kappa = \tau \wedge \phi$. 
Then $\tau (\omega_j) + \phi (\omega_i)=\sigma (\omega_j) + \kappa (\omega_i)$.
\end{lem}
\begin{proof} 
This is just a direct consequence of the fact that 
$\phi + \tau = \phi \vee \tau + \phi \wedge \tau$ in the polytope
described in remark~\ref{rem:polytope}, see~\cite{kn:D-Y}.

It is also a straightforward computation by using the fact that
if $\tau =(t_{1},\ldots ,t_{j})$, then
$$
\tau (\omega_{j})=\omega_{j}-\sum_{k=1}^j (\alpha_k + \cdots + \alpha_{t_k}).
$$
\end{proof}

\begin{lem}\label{lem:redexp1}
Let $j \geq i$ and $\phi \in W_i^w$, $\tau \in W_j^w$ with 
$\sigma = \tau \vee \phi$ and $\kappa = \tau \wedge \phi$. Then, we 
have the following
\begin{enumerate}
\item if $s_{i_1} \cdots s_{i_k} 
\tau = \sigma$ and $\ell (\sigma )=\ell (\tau )+k$, then 
$s_{i_1} \cdots s_{i_k} \breve{\kappa} = \breve{\phi}$ with 
$\ell (\breve{\kappa}) + k = \ell (\breve{\phi})$; or equivalently 
$s_{i_1} \cdots s_{i_k}\kappa = \phi$.
\item if $s_{j_1} \cdots s_{j_l} \breve{\phi} =
\sigma$ with $\ell (\sigma )=\ell (\breve{\phi})+l$, then 
$s_{j_1} \cdots s_{j_l} \breve{\kappa} =\tau$ with
$\ell (\tau )=\ell (\breve {\kappa})+l$.
\item The sets $\{ \alpha_{i_{p}} \}$
and $\{ \alpha_{j_{q}} \}$ has empty intersection and $s_{i_{p}}, 
s_{j_{q}}$ commute.
\end{enumerate}
\end{lem}
\begin{proof} 
Note that as a consequence of Definition~\ref{defn:bar,tilde}, we have 
$\sigma = \tau \vee \breve{\phi}$ and $\breve{\kappa} =
\tau \wedge \breve{\phi}$. Using lemmas~7.17 and~7.18 of 
\cite{kn:G-L}, we conclude that there exist $\alpha_{i_{1}}, \ldots ,
\alpha_{i_{k}}$ and $\alpha_{j_{1}} , \ldots , \alpha_{j_{l}}$ all simple 
enjoying the properties stated above. 
\end{proof}

\section{Theorem on degeneration}
\label{sec:thmdeg}
\setcounter{equation}{0}

Let us recall some basic facts on standard monomials.

Let $\phi\in  W_i$ and $\phi= s_{i_{r}} \cdots s_{i_{1}}$ be a reduced 
expression for $\phi$. Then the vector $Q_{\phi} := X_{-\alpha_{i_{r}}} \cdots 
X_{-\alpha_{i_{1}}} v_{\omega_i}$ is an extremal weight vector in $V(\omega_i)$
of weight $\phi (\omega_i)$. It is shown in \cite{kn:L-S} that
$Q_{\phi}$ is independent of the choice of reduced expression of $\phi$.
Further, we have the following lemmas from~\cite{kn:L-S}:

\begin{lem}\label{lem:Q}
The set $\{ Q_{\tau} \mid \tau \in W_i, \tau \preceq w \}$
is a $\Z$-basis for $E_{\Z, w} (\omega_i)$.
\end{lem}

Let $\{ P_{\tau} \mid \tau \in W_i \}$ be the
$\Z$-basis of $V^{\ast}_{\Z} (\omega_i)$ dual to 
$\{ Q_{\tau} \mid \tau \in W_i \}$. 
Then the set $\{ P_{\tau} \mid \tau \in W_i, \tau \preceq w \}$
is a $\Z$-basis for 
$H^0 (S_{\Z}(w),{\mathcal L}_{\Z, \omega_i})= E^{\ast}_{\Z, w}(\omega_i)$.

\begin{lem}\label{lem:P}
Let $\sigma \succ \kappa \in W_i$ and 
$\sigma = s_{i_{r}} \cdots s_{i_{1}} \kappa$ and 
$\ell (\sigma) - \ell (\kappa) =r$. 
Then we have $P_{\kappa} = (-1)^r X_{-\alpha_{i_{1}}} \cdots X_{-\alpha_{i_{r}}}
P_{\sigma}$.
\end{lem}

For a field $k$, let us denote the canonical image of 
$P_w$ in $H^0 (G/P_i , {\mathcal L}_{\omega_i})$ by $p_w, w \in W_i$.

\begin{defn}\cite{kn:L-S}
A monomial 
$p_{\tau_{r,k_{r}}}\cdots p_{\tau_{r,1}}p_{\tau_{r-1,k_{r-1}}}\cdots 
p_{\tau_{1,1}}$, 
where $\tau_{i,j} \in W_{i}^w$, is called homogeneous of degree 
$(k_1, \ldots , k_r)$
and of total degree $\sum_{j=1}^{r}k_{j}$.

It is called standard on $S(w)$ if for each $i,j$, 
there exists $\tilde{\tau}_{i,j} \in W$, whose class in $W_i$ is $\tau_{i,j}$,
and $\tilde{\tau}_{1,1} \preceq\cdots \preceq \tilde{\tau}_{r,k_{r}}\preceq w$ 
in $W$. In other words $p_{\tau_{r,k_{r}}} \cdots p_{\tau_{1,1}}$
is standard on $S(w)$ if  
$\tau_{1,1} \preceq_w \cdots \preceq_w \tau_{r,k_{r}}\preceq_w w$.
\end{defn}

\begin{thm} \cite{kn:L-S}
\label{thm:standard}
\begin{enumerate}
\item Let $w \in W$. Then, denoting $\overline{w}$
the representative of $w$ in $W_i$, for $\tau \in W_i$,
$p_{\tau} \mid_{S(\overline{w})} \not =0$ if and only if
$\tau \preceq \overline{w}$. Furthermore, $\{ p_{\tau}
\mid \tau \in W_i^w \}$ is a $k$-basis for $H^0 (S(\overline{w}),
{\mathcal L}_{\omega_i})$.

\item The standard monomials on $S(w)$ of degree
$(k_1, \ldots , k_n)$ form a basis of 
$H^0 (S(w), \bigotimes_{i=1}^n {\mathcal L}_{\omega_i}^{\otimes k_i})$.
\end{enumerate}
\end{thm}

Let $H$ be a finite distributive lattice. Denote by
$P = k [x_{\alpha} , \alpha\in H]$ and $I (H) \subset P$ the ideal 
generated by the binomials  
$\{ x_{\alpha} x_{\beta} - x_{\alpha \vee \beta} x_{\alpha \wedge \beta} \mid
\alpha, \beta \in H$ non-comparable$\}$.

Let $R=\bigoplus_{\lambda\, {\rm dominant}}H^{0}(S(w),{\mathcal L}_{\lambda})$ 
be the homogeneous coordinate ring of a multicone over $S_{w}$.
By the previous theorem, $R$ has a basis indexed by standard monomials
on $S(w)$. Thus we have the surjective map $\pi : P \rightarrow R$
sending $x_{\alpha} \mapsto p_{\alpha}$ where we have $H$ is the set of
standard monomials on $S(w)$. Let $I = \ker \pi$ 
which is an ideal generated by relations in total degree $2$ of the form
\begin{equation}
p_{\tau} p_{\phi} - \sum c_{\theta \psi} p_{\theta} p_{\psi}
\end{equation}
where $p_{\tau}p_{\phi}$ is non standard and the $p_{\theta}p_{\psi}$'s 
are standard. These are called straightening relations.

\begin{thm}\cite{kn:G-L}
\label{thm:degeneration}
Assume that ${\mathcal W}^{w}$ is a distributive lattice such that the ideal 
$I$ is generated by the straightening relations of the form 
\begin{equation}\label{eq:straight}
p_{\tau} p_{\phi}-\sum c_{\theta \psi} p_{\theta} p_{\psi}
\end{equation}
where $\tau $, $\phi$ are non-comparable and $\theta \succeq \psi$. 
Further, suppose that we have
\begin{enumerate}
\item $c_{\tau \vee \phi,\tau \wedge \phi}=1$, {\it i.e.} 
$p_{\tau \vee \phi} p_{\tau \wedge \phi}$
occurs on the right-hand side of eq.~$\ref{eq:straight}$
with coefficient $1$.
\item $\tau , \phi \in ]\psi, \theta[ = \{
\gamma \in {\mathcal W}^{w} \mid \psi \preceq \gamma \preceq \theta \}$ for
every pair $(\theta , \psi)$ appearing on the right-hand
side of eq.~$\ref{eq:straight}$.
\item There exist integers $n_1 , \ldots , n_d \geq 1$ and an embedding 
of distributive lattices
$$
\iota :{\mathcal W}^{w} \hookrightarrow  
\bigcup_{d=1}^{n}{\mathcal C} (n_1 , \ldots , n_d)
$$
where ${\mathcal C} (n_1 , \ldots , n_d)$
is the set of $d$-tuples $(i_1, \ldots , i_d)$ with $1 \leq i_j \leq n_j$,
such that for every pair $(\theta , \psi)$ appearing on the 
right-hand side of eq.~$\ref{eq:straight}$,
$\iota(\tau)\dot{\cup}\iota(\phi )=\iota(\theta)\dot{\cup}\iota(\psi )$ 
where $\dot{\cup}$ denotes the disjoint union.
\end{enumerate}
Then there exists a flat deformation whose special fiber is 
$P/I({\mathcal W}^{w})$ and whose general fiber is $R$.
\end{thm}

By lemma \ref{distri}, if $w$ is as in definition~\ref{defn:bar,tilde},
then ${\mathcal W}^{w}$ is a distributive lattice. In the next sections,
we shall prove that all the conditions of the theorem are satisfied.
Let us assume in the next sections that $w$ is as in 
definition~\ref{defn:bar,tilde}.

\section{Condition (2) of Theorem~\ref{thm:degeneration}}
\setcounter{equation}{0}

\begin{thm}\label{thm:strelation}
$($\cite{kn:L-M-S}, \cite{kn:L-S}$)$ 
Let $i\leq j$,  $\tau\in W_{j}^{w}$, $\phi \in W_{i}^{w}$ and 
$p_{\tau} p_{\phi}$ be a non standard monomial on $S(w)$. Let 
the corresponding straightening relation be given by
\begin{equation}
\label{eq:strelation}
p_{\tau} p_{\phi} = \sum_{l=1}^N c_l \ p_{\theta_l} p_{\psi_l}.
\end{equation}
Then $\tau, \phi \prec_w \theta_l$, $\psi_l \prec_w \tau, \phi $ for all $l$
such that $c_l \not = 0$.
\end{thm}

\begin{proof} 
The proof given here is just a generalization of the proof of 
proposition~2.5 of \cite{kn:H-L}.
Among the $\theta_i$ choose a minimal one,
which we denote by $\theta$. Let us reindex the $\theta_l$ so that 
$\theta = \theta_l$ for $1 \leq l \leq s$.
Note that since $\theta$ is minimal we have
$\theta_l \not \preceq \theta$ for $s < l \leq N$.
Since $p_{\theta_l} p_{\psi_l}$ is standard, we can
choose $\kappa^{(l)}_1, \kappa^{(l)}_2 \in W$ such that
$\kappa^{(l)}_2 \preceq \kappa^{(l)}_1\preceq w$, the
class of $\kappa^{(l)}_1$ in $W_j$ is $\theta_l$
and the class of $\kappa^{(l)}_1$ in $W_i$ is $\psi_l$.
Let $Z_1 = \bigcup_{l=1}^s S (\kappa^{(l)}_1)$ and restrict 
eq.~\ref{eq:strelation} to $Z_1$. Then $p_{\theta_l} 
p_{\psi_l} \mid_{Z_1}$ is standard on $Z_1$ for $1 \leq l \leq s$
and $p_{\theta_l} p_{\psi_l} \mid_{Z_1}
\equiv 0$ for $s < l \leq N$. By the linear independence
of standard monomials, eq.~\ref{eq:strelation}
restricted to $Z_1$ is not zero. Hence
$p_{\tau} p_{\phi} \mid_{Z_1} \not = 0$. This implies
that $\tau , \phi \prec \kappa^{(l)}_1$. 
According to Theorem~\ref{stdness} (or
Lemma~8.12 of \cite{kn:D-Y}) we have $\tau , \phi \preceq_w\theta$;
note that $\tau$ (or $\phi$) can not be equal to $\theta$,
because $p_{\tau}p_{\phi}$ is non standard.
>From this argument we deduce that $\tau, \phi \prec_w \theta_l$
for all $l$.

Let $\sigma = \tau \vee \phi \in W_j^w$ and 
$\kappa= \tau \wedge \phi \in W_i^w$. Now $\theta_l \in W_j^w$ and 
$\psi_l \in W^w_i$. By weight consideration, we have 
$\sigma (\omega_j) + \kappa (\omega_i) =\theta_l (\omega_j)+\psi_l (\omega_i)$.
Furthermore $\tau, \phi \prec_w \theta_l$ implies that 
$\sigma \preceq_w \theta_l$, or equivalently $\sigma \preceq \theta_l$ since 
both belong to $W_j$. 
Therefore $\theta_l (\omega_j) \leq \sigma (\omega_j)$, which
implies that $\kappa (\omega_i) \leq \psi_l (\omega_i)$.
Therefore $\psi_l \preceq \kappa$. In other words
$\psi_l \preceq_w \kappa \prec_w \tau, \phi$.
\end{proof}

\begin{cor}
\label{cor:strelation}
Let the notations be as in Lemma~\ref{lem:equalsum1}. Then
in the straightening relation $p_{\tau} p_{\phi} = \sum c_{\theta \psi} 
p_{\theta} p_{\psi}$, either $\sigma \prec_w \theta$ or
$\theta = \sigma$, $\psi = \kappa$.
\end{cor}
\begin{proof} 
>From Theorem~\ref{thm:strelation}, we know
that for any pair $(\theta, \psi)$ on the right-hand side, 
$\sigma\preceq_w \theta$ and $\psi \preceq_w \kappa$. 
Moreover if $\sigma = \theta$, then
due to weight considerations, {\it i.e.} 
$\theta (\omega_j) + \psi (\omega_i) = \sigma (\omega_j) + \kappa (\omega_i)$,
we see that $\kappa = \psi$.
\end{proof}

\section{Condition (3) of Theorem~\ref{thm:degeneration}}
\setcounter{equation}{0}

Considering the set ${\mathcal W}^w := \coprod_{i=1}^n W_i^w$,
we noted at the beginning of Section~\ref{sec:prelim} that
each set $W_i^w$ can be identified with the  subset
of $i$-tuples $(a_1, \ldots , a_i)$ where $0 \leq a_1 <
\cdots < a_i \leq n$ and $(a_1, \ldots , a_i)$ is smaller 
than the representative of $w$ in $W_i$.
Hence we have $\iota : {\mathcal W}^w \hookrightarrow 
\bigcup_{d=1}^n {\mathcal C} (n_1 , \ldots , n_d)$. 
For simplicity, we shall denote $\iota (\tau )$ also by $\tau$.
We want to prove the following lemma:

\begin{lem}\label{weightcondition}
Let $\tau, \phi$ be two non-comparable elements in
${\mathcal W}^w$. Then for any $(\theta, \psi )$ appearing
on the right-hand side of the straightening relation~\ref{eq:straight}, 
$\theta \dot{\cup} \psi = \tau \dot{\cup} \phi$. 
\end{lem}

\begin{proof}
Let $\tau = (t_1, \ldots , t_j)$, $\phi =
(r_1, \ldots , r_i)$, $ \theta = (a_1, \ldots , a_j)$
and $\psi = (b_1, \ldots ,b_i)$. A
necessary condition for $p_{\theta} p_{\psi}$ to
appear on the right-hand side of eq.~\ref{eq:straight} is
$\tau (\omega_j) + \phi (\omega_i) =\theta (\omega_j) + \psi (\omega_i)$.
Here, we shall prove that this condition
immediately implies the assertion, {\it i.e.} $\{ t_1, \ldots ,
t_j \} \dot{\cup} \{r_1, \ldots , r_i \} = \{ a_1, \ldots , a_j,
\} \dot{\cup} \{b_1, \ldots , b_i \}$. The proof is by induction on
$i+j$. 

The fact that $\tau (\omega_j) + \phi (\omega_i) =
\theta (\omega_j) + \psi (\omega_i)$ implies (see the proof of 
lemma~\ref{lem:equalsum1})
\begin{equation}\label{eq:weight=}
\begin{array}{l}
\displaystyle\sum_{k=1}^j (\alpha_k + \cdots + \alpha_{t_k}) +
\sum_{l=1}^i (\alpha_l + \cdots + \alpha_{r_l}) =\\
\hskip5em \displaystyle\sum_{k=1}^j (\alpha_k + \cdots + \alpha_{a_k}) +
\sum_{l=1}^i (\alpha_l + \cdots + \alpha_{b_l})
\end{array}
\end{equation}
Note that 
$$
\max \{ t_1, \ldots ,t_j, r_1, \ldots , r_i \} = \max \{ t_j, r_i \}
$$ 
and that 
$$
\max \{ a_1, \ldots , a_j, b_1, \ldots , b_i \} = \max \{ a_j, b_i \}.
$$
Then due to the equality in eq.~\ref{eq:weight=},
we must have $\max \{ t_j, r_i \} = \max \{ a_j, b_i \} $.
There are four cases to consider.

$\bullet$ Case (1) $t_j = a_j \geq b_i$. This implies that 
$\alpha_j + \cdots + \alpha_{t_j} =\alpha_j + \cdots + \alpha_{a_j}$.
Hence denoting $\tau^{\prime} =(t_1 ,\ldots , t_{j-1})$ and 
$\theta^{\prime} =(a_1, \ldots , a_{j-1})$, 
eq.~\ref{eq:weight=} implies that
$\tau^{\prime} (\omega_{j-1}) + \phi (\omega_i) =
\theta^{\prime} (\omega_{j-1}) + \psi (\omega_i)$.
By induction we are done.

$\bullet$ Case (2) $t_j = b_i > a_j$. Let $m$ be the smallest
number such that $a_{j-m} > b_{i-m}$ (if such an $m$ less than
$i-1$ does not exist, let $m=i$). Note that $b_{i-m+1} \geq
a_{j-m+1} > a_{j-m} > b_{i-m}$. 
Set
\begin{itemize}
\item[] $\tau^{\prime}  = (t_1, \ldots , t_{j-1}) \in W_{j-1}$, 
\item[] $\theta^{\prime} = (a_1, \ldots , a_{j-m}, b_{i-m+1}, b_{i-m+2},
\ldots , b_{i-1}) \in W_{j-1}$,
\item[] $\psi^{\prime} = (b_1,\ldots ,b_{i-m},a_{j-m+1},a_{j-m+2}, \ldots , a_j)
\in W_i$ if $m\neq i$ and 
\item[] $\psi^{\prime} = (a_{j-i+1}, \ldots ,a_j)$ if $m=i$. 
\end{itemize}
Since $i \leq j$, we have $i-k-1 \leq j-k-1 \leq a_{j-k}$ for $0 \leq k < i$.
Therefore $\psi'\in W_i$. Using the fact that
for $0 \leq k < m$, we have $i-k-1 \leq j-k-1 \leq a_{j-k} \leq b_{i-k}$, then
\begin{equation}
\begin{array}{l}\label{eq:sum1}
(\alpha_{j-k} + \cdots + \alpha_{a_{j-k}}) +
(\alpha_{i-k} + \cdots + \alpha_{b_{i-k}}) =\\ 
\hskip2em (\alpha_{j-k} + \cdots + \alpha_{a_{j-k}} +
\alpha_{a_{j-k}+1} + \cdots + \alpha_{b_{i-k}}) +  
(\alpha_{i-k} + \cdots + \alpha_{a_{j-k}})
\end{array}
\end{equation}
>From eqs.~\ref{eq:weight=} and~\ref{eq:sum1}, we can conclude
that $\tau^{\prime} (\omega_{j-1}) + \phi (\omega_i) =
\theta^{\prime} (\omega_{j-1}) + \psi^{\prime} (\omega_i)$.
The rest follows by induction. 

$\bullet$ Case (3) $r_i = b_i \geq a_j$ is similar to case (1). 

$\bullet$ Case (4) $r_i = a_j > b_i$ is similar to case (2).
\end{proof}

In fact, we have proved:

\begin{lem}
Let $j\geq i$, $\tau ,\theta\in W_{j}$, $\phi ,\psi\in W_{i}$ be
such that $\tau (\omega_j) +  \phi (\omega_{i}) = \theta 
(\omega_j) +  \psi (\omega_{i})$. Then 
$\theta \dot{\cup} \psi = \tau \dot{\cup} \phi$. 
\end{lem}

\section{Condition (1) of Theorem~\ref{thm:degeneration}}

\begin{prop}
\label{prop:+-1}
Let $\tau, \phi \in \mathcal{W}^w$ be two non-comparable elements. 
Then in the straightening relation~$\ref{eq:strelation}$,
$p_{\tau \vee \phi} p_{\tau \wedge \phi}$ occurs with
coefficient $\pm 1$.
\end{prop}
\begin{proof} 
As before, denote $\sigma = \tau \vee \phi$,
$\kappa = \tau \wedge \phi$. Note that $\tau, \phi \prec_w
\sigma$ (that is there exist liftings $\tilde{\tau}, 
\tilde{\phi}, \tilde{\sigma}$ in $W$ such that
$\tilde{\tau}, \tilde{\phi} \preceq \tilde{\sigma}\preceq w$). 
Corollary~\ref{cor:strelation} implies 
that the restriction of eq.~\ref{eq:strelation}
to the Schubert variety $S (\tilde{\sigma})$ is 
$p_{\tau} p_{\phi} = a p_{\sigma} p_{\kappa}$, where
$a \not = 0$. Since standard monomial basis is characteristic free, 
this holds in any characteristics. Hence $a = \pm 1$.
\end{proof}

So now we have to prove that $a =1$. 
Since the irreducible representation $V(\omega_i + \omega_j)$,
appears as a direct sum in the decomposition 
in $V(\omega_j) \otimes V(\omega_i)$ into irreducible 
representations, we have an imbedding $V(\omega_i + \omega_j)
\hookrightarrow V(\omega_j) \otimes V(\omega_i)$. Note that
since the weight space of weight $\omega_i + \omega_j$ is one-dimensional,
the element $v_{\omega_j} \otimes v_{\omega_i}$
belongs to $V(\omega_i + \omega_j)$. The imbedding above
induces a projection $H^0 (G/B, {\mathcal L}_{\omega_i})
\otimes H^0 (G/B, {\mathcal L}_{\omega_j}) \rightarrow
H^0 (G/B, {\mathcal L}_{\omega_i} \otimes {\mathcal L}_{\omega_j})$.
For simplicity we shall denote the image of $f \otimes g$ under
this projection by $fg$. We shall construct a basis for
$E_{\Z,w}(\omega_{i}+\omega_{j})$ which is a ``rank two'' version
of the one given in \cite{kn:L-S}.

In the following let $i \leq j$ (that is no element of 
$W_i^w$ can be bigger than an element of $W_j^w$) and recall  
from lemma~\ref{lem:Q} that, for $\phi \in W_i$, we have denoted by
$Q_{\phi} $ an extremal weight vector in  $V_{\Z}(\omega_i)$ of weight
$\phi (\omega_i)$.

Let $\Sigma(w):=\{ (\tau,\sigma )\in W_{j}^{w}\times W_{i}^{w}$ $\mid$
there exist liftings $\widetilde{\tau},\widetilde{\sigma}$ in $W$ such that 
$\widetilde{\sigma}\preceq\widetilde{\tau}\preceq w \}$.   

\begin{defn}
\label{defn:e-s,k}
Let $w$ be as in Definition~\ref{defn:bar,tilde}.
Let $\kappa \in W_i^w$, $\sigma \in W_j^w$ be such that
$(\sigma ,\kappa )\in\Sigma (w)$ and let $\sigma = s_{i_{r}}
\cdots s_{i_{1}} \breve{\kappa}$ where 
$r = \ell (\sigma) - \ell (\breve{\kappa})$. Define $E_{\breve{\kappa},
\kappa} := Q_{\breve{\kappa}} \otimes Q_{\kappa} \in
V_{\Z}(\omega_j) \otimes V_{\Z}(\omega_i)$ and
define $E_{\sigma , \kappa} := X_{-\alpha_{i_{r}}} \cdots X_{-\alpha_{i_{1}}}
E_{\breve{\kappa}, \kappa}$.
\end{defn}

Note that $E_{\breve{\kappa},\kappa}$ is an extremal weight vector since
$\breve{\kappa}$ is the image of $\kappa$ (the minimal representative in $W$) 
in $W_{j}$. 
It is also clear that $E_{\sigma , \kappa}$ is a weight vector of weight
$\kappa (\omega_i) + \sigma (\omega_j)$.

\begin{prop}\label{pro:basis}
Let $w\in W$ be as in definition \ref{defn:bar,tilde}. Then
$E_{\sigma ,\kappa}$ does not depend on the choice of reduced expression and
the set $\{ E_{\sigma , \kappa} \mid 
\kappa \in W_i^{w}, \sigma \in W_j^{w}, \kappa \preceq_w \sigma\}$
is a $\Z$-basis for the Demazure module
$E_{\Z, w} (\omega_i + \omega_j)$.
\end{prop}
\begin{proof}
Let $\sigma = s_{i_{r}}
\cdots s_{i_{1}} \breve{\kappa}=s_{j_{r}}
\cdots s_{j_{1}} \breve{\kappa}$. 
Denote by $\phi = s_{j_{r-1}}\cdots s_{j_{1}}\breve{\kappa}$.
Then we have $\sigma = s_{j_{r}}\phi$. 
Now if $i_{r}=j_{r}$, then we proceed by induction on the length
of $\sigma$. Otherwise, let $k$ be the largest integer such that $j_{r}=i_{k}$.
We have $\phi \vee \breve{\kappa}=\sigma$, thus
by lemma~\ref{lem:redexp1}, we have that $s_{j_{r}}$ commute
with $s_{i_{l}}$ for $l\geq k$. 
Thus 
$$
X_{-\alpha_{i_{r}}}\cdots X_{-\alpha_{1}}E_{\breve{\kappa},\kappa}
=X_{-\alpha_{i_{k}}}X_{-\alpha_{i_{r}}}\cdots 
X_{-\alpha_{i_{k+1}}}X_{-\alpha_{i_{k-1}}}\cdots X_{-\alpha_{i_{1}}}
E_{\breve{\kappa},\kappa}
$$
By induction, 
$E_{\phi ,\kappa}=X_{-\alpha_{j_{r-1}}}\cdots X_{-\alpha_{j_{1}}}
E_{\breve{\kappa},\kappa}$.
Therefore the right-hand side is $X_{-\alpha_{j_{r}}}E_{\phi ,\kappa}$
and we have proved that the definition of
$E_{\sigma ,\kappa}$ does not depend on the choice of the reduced expression. 

We are left to show that these elements form a basis for 
$E_{\Z,w}(\omega_{i}+\omega_{j})$.

We claim that $E_{\sigma,\kappa}\in E_{\Z,w}(\omega_{i}+\omega_{j})$.
It is clear that $E_{\breve{\kappa},\kappa}\in 
E_{\Z,\kappa}(\omega_{i}+\omega_{j})$. Now, since $w$ satisfies the
condition of definition \ref{defn:bar,tilde}, we have 
$w\succeq s_{i_{r}}\cdots s_{i_{1}}\kappa$, thus
$$
E_{\sigma,\kappa}\in X_{-\alpha_{i_{r}}}\cdots X_{-\alpha_{i_{1}}}
E_{\Z,\kappa}(\omega_{i}+\omega_{j})\subset E_{\Z,w}(\omega_{i}+\omega_{j}).
$$
We have therefore our claim.

Now by the definition of $E_{\sigma,\kappa}$, we have
$$
E_{\sigma,\kappa}=Q_{\sigma}\otimes Q_{\kappa}+
\sum_{(u,v)\in I} Q_{u}\otimes Q_{v}
$$
where $I\subset W_{j}\times W_{i}$ and for each $(u,v)\in I$, we have
$u \prec \sigma $ in $W_{j}$, $v \succ \kappa$ in $W_{i}$
and $\sigma (\omega_{j})+\kappa (\omega_{i}) = u(\omega_{j})+v(\omega_{i})$.
It is now clear that the $E_{\sigma,\kappa}$'s are independant.

Further, one deduces from the expression for $E_{\sigma,\kappa}$ above
that the $\Z$-submodule generated by the 
$E_{\sigma,\kappa}$'s is a direct summand of the tensor product
$V_{\Z}(\omega_{j})\otimes V_{\Z}(\omega_{i})$.
Finally, by standard monomial theory, the cardinal of $\Sigma (w)$ 
is the rank of $E_{\Z,w}(\omega_{i}+\omega_{j})$.
So the result follows.
\end{proof}

We can now prove that $a=1$.

\begin{cor}\label{cor:a=1}
Let the notations be as in Lemma~\ref{lem:equalsum1},
then in the straightening relation 
$p_{\tau} p_{\phi} = \sum_{l=1}^N c_l \
p_{\theta_l} p_{\psi_l}$, the term $p_{\sigma} p_{\kappa}$
occurs on the right hand side with coefficient $1$.
\end{cor}
\begin{proof}
Recall from the proof of~\ref{pro:basis} that 
$$
E_{\sigma,\kappa}=Q_{\sigma}\otimes Q_{\kappa}+
\sum_{(u,v)\in I} Q_{u}\otimes Q_{v}
$$
where $I\subset W_{j}\times W_{i}$ and for each $(u,v)\in I$, we have
$u \prec \sigma $ in $W_{j}$, $v \succ \kappa$ in $W_{i}$
and $\sigma (\omega_{j})+\kappa (\omega_{i}) = u(\omega_{j})+v(\omega_{i})$.

Let us apply $p_{\tau}p_{\phi}$ 
to $E_{\sigma,\kappa}$. Then from the explicite expression of
$E_{\sigma,\kappa}$ above, this is either $0$ or $1$
depending if $Q_{\tau}\times Q_{\phi}$ appears in the right hand
side or not.

On the other hand, if we replace $p_{\tau}p_{\phi}$ by
the right hand side of the straightening relation, then it is clear from
Theorem~\ref{thm:strelation} that
the same evaluation yields $a_{\sigma,\kappa}$. But this
is non zero from Proposition~\ref{prop:+-1}. So it must be $1$.
\end{proof}

\section{Consequence}
\label{consequence}

As an immediate consequence, we have:
\begin{thm}
Let $w$ be as in definition~\ref{defn:bar,tilde}. Then
there exists a flat deformation whose special fiber is a toric 
variety and whose general fiber is $S(w)$.
\end{thm}
\begin{proof}
By theorem~\ref{thm:standard}, there exists a flat deformation
whose general fiber is $S(w)$ and whose special fiber is a variety defined 
by a binomial ideal associated to a distributive lattice.
This latter is toric as shown in \cite{kn:E-S}.
\end{proof}

\begin{rem}
If we look closely at the proofs, then we realize that 
theorem~\ref{thm:standard} can be replaced by the following.

Suppose that ${\mathcal W}^{w}$ admits a structure of distributive
lattice such that 
\begin{enumerate}
\item the partial order corresponds to standardness, cf theorem~\ref{stdness};
\item weights are preserved, cf lemma~\ref{lem:equalsum1};
\item lemma~\ref{lem:redexp1} is satisfied.
\end{enumerate}
Then there exists a flat deformation whose special fiber is a toric 
variety and whose general fiber is $S(w)$.
\end{rem}

In particular, consider the bijection $\Theta$ of $W$ defined by
$s_{i}\mapsto s_{n+1-i}$ induced by the non trivial Dynkin diagram
automorphism. This induces a bijection
between $W_{i}$ and $W_{n+1-i}$ which preverses the induced Bruhat
order. Now let $w$ be as in Definition~\ref{defn:bar,tilde}, then
$\Theta$ induces a structure of distributive lattice
on ${\mathcal W}^{\Theta(w)}$.
It is easy to check that the same proof works. Thus we
have,

\begin{thm}\label{lastthm}
Let $w$ or $\Theta(w)$ be as in definition~\ref{defn:bar,tilde}. Then
there exists a flat deformation whose special fiber is a toric 
variety and whose general fiber is $S(w)$.
\end{thm}

\begin{rem}\label{connectedparts}
As noticed in~\cite{kn:D-Y}, we can extend our results to the following
elements. Let $0\leq k_{1}< k_{2}<\cdots < k_{r+1}\leq n+1$, and for 
$1\leq i\leq r$, let $S_{i}$ be the subgroup of $W$ generated by the 
reflections $s_{k_{i}+1},\cdots s_{k_{i+1}-1}$. 

Now suppose that $w=w_{1}\cdots w_{r}$ where $w_{i}\in S_{i}$. Then
it is clear that $w_{i}$ and $w_{j}$ commute if $i\neq j$ and
it follows easily that if each $w_{i}$ satisfies the condition 
of theorem~\ref{lastthm}, {\it i.e.} either $w_{i}$ or $\Theta (w_{i})$
is as in definition~\ref{defn:bar,tilde}, then the conclusion of the same 
theorem holds for $w$.

For example, the element $s_{1}s_{2}s_{5}s_{4}$ satisfies the above 
conditions.
\end{rem}

Our results apply to all the elements of $W$ in the case of $SL_{3}$
thus giving a more general proof to~\cite{kn:D}.
However, in the case of $SL_{4}$, there are precisely $4$ elements
for which the condition of the theorem is not satisfied. 
Namely, they are $s_{1}s_{3}s_{2}$, $s_{2}s_{1}s_{3}$, $s_{2}s_{1}s_{3}s_{2}$
and $s_{1}s_{2}s_{3}s_{2}s_{1}$. The main problem in these cases is
that standardness is not transitive in all the obvious ``orderings''.

\end{document}